\documentclass[reqno]{amsart}%
\usepackage[english]{babel}
\usepackage{amscd}
\usepackage{amsmath}
\usepackage{graphicx}
\usepackage{amsfonts}
\usepackage{amssymb}%

\makeatletter \DeclareMathSymbol{\Gamma}{\mathalpha}{letters}{"00}
\DeclareMathSymbol{\Delta}{\mathalpha}{letters}{"01}
\DeclareMathSymbol{\Theta}{\mathalpha}{letters}{"02}
\DeclareMathSymbol{\Lambda}{\mathalpha}{letters}{"03}
\DeclareMathSymbol{\Xi}{\mathalpha}{letters}{"04}
\DeclareMathSymbol{\Pi}{\mathalpha}{letters}{"05}
\DeclareMathSymbol{\Sigma}{\mathalpha}{letters}{"06}
\DeclareMathSymbol{\Upsilon}{\mathalpha}{letters}{"07}
\DeclareMathSymbol{\Phi}{\mathalpha}{letters}{"08}
\DeclareMathSymbol{\Psi}{\mathalpha}{letters}{"09}
\DeclareMathSymbol{\Omega}{\mathalpha}{letters}{"0A}
\DeclareMathSymbol{\varGamma}{\mathalpha}{operators}{"00}
\DeclareMathSymbol{\varDelta}{\mathalpha}{operators}{"01}
\DeclareMathSymbol{\varTheta}{\mathalpha}{operators}{"02}
\DeclareMathSymbol{\varLambda}{\mathalpha}{operators}{"03}
\DeclareMathSymbol{\varXi}{\mathalpha}{operators}{"04}
\DeclareMathSymbol{\varPi}{\mathalpha}{operators}{"05}
\DeclareMathSymbol{\varSigma}{\mathalpha}{operators}{"06}
\DeclareMathSymbol{\varUpsilon}{\mathalpha}{operators}{"07}
\DeclareMathSymbol{\varPhi}{\mathalpha}{operators}{"08}
\DeclareMathSymbol{\varPsi}{\mathalpha}{operators}{"09}
\DeclareMathSymbol{\varOmega}{\mathalpha}{operators}{"0A}
\newcommand{\allmodesymb}[2]{\relax\ifmmode{\mathchoice
{\mbox{\fontsize{\tf@size}{\tf@size}#1{#2}}}
{\mbox{\fontsize{\tf@size}{\tf@size}#1{#2}}}
{\mbox{\fontsize{\sf@size}{\sf@size}#1{#2}}}
{\mbox{\fontsize{\ssf@size}{\ssf@size}#1{#2}}}} \else
\mbox{#1{#2}}\fi}
\makeatother
\makeatletter
\renewcommand*\subjclass[2][2000]{%
  \def\@subjclass{#2}%
  \@ifundefined{subjclassname@#1}{%
    \ClassWarning{\@classname}{Unknown edition (#1) of Mathematics%
      Subject Classification; using '2000'.}%
  }{%
    \@xp\let\@xp\subjclassname\csname subjclassname@#1\endcsname%
  }%
} \makeatother

\theoremstyle{plain}
\newtheorem{theorem}{Theorem}[section]
\newtheorem{corollary}[theorem]{Corollary}

\newtheorem{proposition}[theorem]{Proposition}
\newtheorem{definition}{Definition}[section]

\theoremstyle{remark}
\newtheorem*{remark}{Remark}

\numberwithin{equation}{section} \allowdisplaybreaks
\begin{document}
\author{Guy Roos}
\address{Nevski prospekt 113/4-53, 191024 St Petersburg, Russian Federation}
\email{guy.roos@normalesup.org}
\thanks{One hour lecture for graduate students, SCV2004, Beijing}
\title{Weighted Bergman kernels and virtual Bergman kernels}
\date{August 25, 2004}

\begin{abstract}
We introduce the notion of \textit{virtual Bergman kernel} and study some of
its applications.

\end{abstract}
\maketitle

\tableofcontents

\section*{Introduction}

Let $\Omega\subset\mathbb{C}^{n}$ be a domain and $p:\Omega\rightarrow
]0,+\infty\lbrack$ a weight function on $\Omega$. Consider the
\textquotedblleft inflated domains\textquotedblright\
\begin{align}
\widehat{\Omega}_{1}  &  =\left\{  \left(  z,\zeta\right)  \in\Omega
\times\mathbb{C}\mid\left\vert \zeta\right\vert ^{2}<p(z)\right\}  ,\\
\widehat{\Omega}_{m}  &  =\left\{  \left(  z,Z\right)  \in\Omega
\times\mathbb{C}^{m}\mid\left\Vert Z\right\Vert ^{2}<p(z)\right\}  ,
\end{align}
where $\left\Vert ~\right\Vert $ is the standard Hermitian norm on
$\mathbb{C}^{m}$.

In our joint work \cite{YinRoos2003} with Yin Weiping, we computed explicitly
the Bergman kernel of some \textquotedblleft egg domains\textquotedblright;
among them $\widehat{\Omega}_{1}$, when $\Omega$ is a bounded symmetric domain
and $p$ a real power of the generic norm of $\Omega$. We then obtained the
Bergman kernel of the corresponding $\widehat{\Omega}_{m}$ by using the
\emph{\textquotedblleft inflation principle\textquotedblright} of
\cite{BoasFuStraube1999}, which allows to deduce (for any weight function $p$)
the Bergman kernel of $\widehat{\Omega}_{m}$ from the Bergman kernel of
$\widehat{\Omega}_{1}$. The \textquotedblleft inflation
principle\textquotedblright\ says that if the Bergman kernel of $\widehat
{\Omega}_{1}$ is {%
\[
\widehat{\mathcal{K}}_{1}(z,\zeta)=\mathcal{L}_{1}\left(  z,\left\vert
\zeta\right\vert ^{2}\right)  ,
\]
}then the Bergman kernel of $\widehat{\Omega}_{m}$ is
\begin{equation}
\widehat{\mathcal{K}}_{m}(z,Z)=\frac{1}{m!}\left.  \frac{\partial^{m-1}%
}{\partial r^{m-1}}\mathcal{L}_{1}(z,r)\right\vert _{r=\left\Vert Z\right\Vert
^{2}}. \label{VB0}%
\end{equation}

It appears that the two previous steps can be unified in the following way.
There exists a function $\mathcal{L}_{0}(z,r)$, defined in a neighborhood of
$\Omega\times\{0\}$ in $\Omega\times\lbrack0,+\infty\lbrack$, such that for
all $m\geq1$, the Bergman kernel of $\widehat{\Omega}_{m}$ is
\begin{equation}
\widehat{\mathcal{K}}_{m}(z,Z)=\frac{1}{m!}\left.  \frac{\partial^{m}%
}{\partial r^{m}}\mathcal{L}_{0}(z,r)\right\vert _{r=\left\Vert Z\right\Vert
^{2}}.
\end{equation}
We call $\mathcal{L}_{0}(z,r)$ the \emph{virtual Bergman kernel }of $\left(
\Omega,p\right)  $. Its existence is closely related to the \textquotedblleft
Forelli-Rudin construction\textquotedblright\ (\cite{ForelliRudin1974},
\cite{Ligocka1989}, \cite{Englis2000}).

In this talk, we investigate the properties of this virtual Bergman kernel. We
then show how it can be explicitly computed on bounded symmetric domains, for
a special but natural choice of the weight function $p$: $p$ is taken to be a
real power of the \textquotedblleft generic norm\textquotedblright\ of the
bounded symmetric domain. The explicit computation of the virtual Bergman
kernel is then related to properties of the Hua integral.

\section{Virtual Bergman kernels}

\subsection{Notations}

Let $V\cong\mathbb{C}^{n}$ be a Hermitian vector space, with Hermitian norm
$\left\Vert ~\right\Vert _{V}$ and volume form $\omega_{V}(z)=\left(
\frac{\operatorname*{i}}{2\pi}\partial\overline{\partial}\left\Vert
z\right\Vert ^{2}\right)  ^{n}$. Let $\Omega$ be a domain in $V$ and
$p:\Omega\rightarrow]0,+\infty\lbrack$ a continuous function on $\Omega$. The
space of holomorphic functions on $\Omega$ is denoted by $\operatorname{Hol}%
(\Omega)$. We denote by $H(\Omega)$ the Bergman space{%
\[
H(\Omega)=H\left(  \Omega,\omega_{V}\right)  =\left\{  f\in\operatorname{Hol}%
(\Omega)\mid\left\Vert f\right\Vert _{\Omega}^{2}=\int_{\Omega}\left\vert
f(z)\right\vert ^{2}\omega_{V}(z)<\infty\right\}
\]
and by }$H(\Omega,p)$ the weighted Bergman space{%
\[
H(\Omega,p)=H\left(  \Omega,p\omega_{V}\right)  =\left\{  f\in
\operatorname{Hol}(\Omega)\mid\left\Vert f\right\Vert _{\Omega,p}^{2}%
=\int_{\Omega}\left\vert f(z)\right\vert ^{2}p(z)\omega_{V}(z)<\infty\right\}
.
\]
The Hilbert products on these spaces are denoted respectively }$\left(
~\mid~\right)  _{\Omega}$ and $\left(  ~\mid~\right)  _{\Omega,p}$. The
Bergman kernel of $\Omega$ (reproducing kernel of $H(\Omega)$) is denoted by
$K_{\Omega}(z,t)$; it is fully determined by $\mathcal{K}_{\Omega}$:%
\[
\mathcal{K}_{\Omega}(z)=K_{\Omega}(z,z)\qquad\left(  z\in\Omega\right)  ,
\]
which we call also Bergman kernel of $\Omega$. In the same way, the weighted
Bergman kernel of $\left(  \Omega,p\right)  $ (reproducing kernel of
$H(\Omega,p)$) is denoted by $K_{\Omega,p}(z,t)$ and is determined by
$\mathcal{K}_{\Omega,p}(z)=K_{\Omega,p}(z,z)$.

\subsection{}

Let
\[
\widehat{\Omega}_{1}=\left\{  \left(  z,\zeta\right)  \in\Omega\times
\mathbb{C}\mid\left|  \zeta\right|  ^{2}<p(z)\right\}  .
\]
We endow $\widehat{\Omega}_{1}$ with the volume form%
\[
\omega_{V}(z)\wedge\omega_{1}(\zeta),
\]
where $\omega_{1}(\zeta)=\frac{\operatorname*{i}}{2\pi}\partial\overline
{\partial}\left|  z\right|  ^{2}$. The Bergman space $H\left(  \widehat
{\Omega}_{1}\right)  $ will be defined with respect to this volume form.

Consider a holomorphic function $f\in\operatorname{Hol}\left(  \widehat
{\Omega}_{1}\right)  $; such a function may be written
\[
f(z,\zeta)=\sum_{k=0}^{\infty}f_{k}(z)\zeta^{k},
\]
with $f_{k}\in\operatorname{Hol}\Omega$. We compute $\left\|  f\right\|
_{\widehat{\Omega}_{1}}^{2}$:%
\begin{align*}
\left\|  f\right\|  _{\widehat{\Omega}_{1}}^{2}  &  =\int_{\widehat{\Omega
}_{1}}\left|  f(z,\zeta)\right|  ^{2}\omega_{V}(z)\wedge\omega_{1}(\zeta)\\
&  =\int_{\Omega}\omega_{V}(z)\left(  \int_{\left|  \zeta\right|  ^{2}%
<p(z)}\left|  f(z,\zeta)\right|  ^{2}\omega_{1}(\zeta)\right) \\
&  =\int_{\Omega}\omega_{V}(z)\left(  \sum_{k=0}^{\infty}\left|
f_{k}(z)\right|  ^{2}\frac{p^{k+1}(z)}{k+1}\right)  ,
\end{align*}
which gives%
\begin{equation}
\left\|  f\right\|  _{\widehat{\Omega}_{1}}^{2}=\sum_{k=0}^{\infty}\frac
{1}{k+1}\left\|  f_{k}\right\|  _{\Omega,p^{k+1}}^{2}. \label{VB1}%
\end{equation}

\begin{proposition}
Let $\widehat{\mathcal{K}}_{1}$ be the Bergman kernel of $\widehat{\Omega}%
_{1}$ and $\mathcal{K}_{\Omega,p^{k}}$ the weighted Bergman kernel of $\Omega$
for the weight function $p^{k}$. Then
\begin{equation}
\widehat{\mathcal{K}}_{1}(z,\zeta)=\sum_{k=0}^{\infty}\left(  k+1\right)
\mathcal{K}_{\Omega,p^{k+1}}(z)\left\vert \zeta\right\vert ^{2k}. \label{VB2}%
\end{equation}

\end{proposition}

\proof For each $k\in\mathbb{N}$, let $\left(  \varphi_{jk}\right)  _{j\in
J_{k}}$ be a Hilbert basis (complete orthonormal system) of $H\left(
\Omega,p^{k}\right)  $. Then it follows from (\ref{VB1}) that
\[
\left(  \left(  k+1\right)  ^{1/2}\varphi_{j,k+1}(z)\zeta^{k}\right)
_{k\in\mathbb{N},\ j\in J_{k+1}}%
\]
is a Hilbert basis of $H\left(  \widehat{\Omega}_{1}\right)  $. From the
classical properties of Bergman kernels, we get
\begin{align*}
\mathcal{K}_{\Omega,p^{k}}(z)  &  =\sum_{j\in J_{k}}\left\vert \varphi
_{jk}(z)\right\vert ^{2},\\
\widehat{\mathcal{K}}_{1}(z,\zeta)  &  =\sum_{k\in\mathbb{N},\,j\in J_{k}%
}\left(  k+1\right)  \left\vert \varphi_{j,k+1}(z)\right\vert ^{2}\left\vert
\zeta\right\vert ^{2k}\\
&  =\sum_{k=0}^{\infty}\left(  k+1\right)  \mathcal{K}_{\Omega,p^{k+1}%
}(z)\left\vert \zeta\right\vert ^{2k}.
\end{align*}
\endproof

\subsection{}

This leads to the following definition.

\begin{definition}
Let $\Omega$ be a domain in $V$ and $p:\Omega\rightarrow]0,+\infty\lbrack$ a
continuous function on $\Omega$. Denote by $K_{\Omega,p^{k}}(z,w)$
($\mathcal{K}_{\Omega,p^{k}}(z)$) the weighted Bergman kernel of $\Omega$ w.r.
to $p^{k}$. The \emph{virtual Bergman kernel} of $\left(  \Omega,p\right)  $
is defined by
\begin{equation}
L_{\Omega,p}\left(  z,w;r\right)  =L_{0}\left(  z,w;r\right)  =\sum
_{k=0}^{\infty}K_{\Omega,p^{k}}(z,w)r^{k}. \label{VB3'}%
\end{equation}
The function $\mathcal{L}_{0}(z,r)=L_{0}(z,z;r)$, i.e.%
\begin{equation}
\mathcal{L}_{\Omega,p}\left(  z,r\right)  =\mathcal{L}_{0}\left(  z,r\right)
=\sum_{k=0}^{\infty}\mathcal{K}_{\Omega,p^{k}}(z)r^{k} \label{VB3}%
\end{equation}
will also be called virtual Bergman kernel of $\left(  \Omega,p\right)  $.
\end{definition}

With these definitions, the relation (\ref{VB2}) may be rewritten and the
reproducing kernel of $\widehat{\Omega}_{1}$ is given by
\begin{align}
\widehat{K}_{1}\left(  \left(  z,\zeta\right)  ,\left(  w,\eta\right)
\right)   &  =L_{1}\left(  z,w;\zeta\overline{\eta}\right)  ,\label{VB4'}\\
\widehat{\mathcal{K}}_{1}(z,\zeta)  &  =\mathcal{L}_{1}\left(  z,\left\vert
\zeta\right\vert ^{2}\right)  , \label{VB4}%
\end{align}
with%
\begin{align*}
L_{1}(z,w;r)  &  =\frac{\partial}{\partial r}L_{0}\left(  z,w;r\right)  ,\\
\mathcal{L}_{1}(z,r)  &  =\frac{\partial}{\partial r}\mathcal{L}_{0}\left(
z,r\right)  .
\end{align*}

\begin{remark}
From the virtual Bergman kernel of $\left(  \Omega,p\right)  $, it is easy to
recover the weighted Bergman kernels of $\left(  \Omega,p^{k}\right)  $
($k\in\mathbb{N}$):%
\begin{equation}
\mathcal{K}_{\Omega,p^{k}}(z)=\frac{1}{k!}\left.  \frac{\partial^{k}}{\partial
r^{k}}\mathcal{L}_{\Omega,p}\left(  z,r\right)  \right\vert _{r=0}.
\label{VB13}%
\end{equation}

\end{remark}

\subsection{}

Let us recall some facts about harmonic analysis in the Hermitian unit ball
$B_{m}$. Let $H\left(  B_{m}\right)  =H\left(  B_{m},\omega_{m}\right)  $ be
the Bergman space of $B_{m}$; let $K_{B_{m}}(Z,T)$, $\mathcal{K}_{B_{m}%
}(Z)=K_{B_{m}}(Z,Z)$ be the Bergman kernel. It is well known (using for
instance the automorphisms of $B_{m}$) that
\begin{equation}
\mathcal{K}_{B_{m}}(Z)=\frac{1}{\left(  1-\left\|  Z\right\|  ^{2}\right)
^{m+1}}. \label{VB6}%
\end{equation}
This may also be written%
\begin{equation}
\mathcal{K}_{B_{m}}(Z)=\frac{1}{m!}\left.  \frac{\partial^{m}}{\partial r^{m}%
}\left(  \frac{1}{1-r}\right)  \right|  _{r=\left\|  Z\right\|  ^{2}}.
\label{VB8}%
\end{equation}

Let $f\in H\left(  B_{m}\right)  $; then $f$ may be written%
\begin{equation}
f(Z)=\sum_{k=0}^{\infty}f_{k}(Z), \label{VB11}%
\end{equation}
where the $f_{k}$ are $k$-homogeneous polynomials, which can be obtained
through%
\[
f_{k}(Z)=\int_{0}^{1}f\left(  \operatorname*{e}\nolimits^{2\pi
\operatorname*{i}\theta}Z\right)  \operatorname*{e}\nolimits^{-2\pi
\operatorname*{i}k\theta}\operatorname{d}\theta.
\]
The expansion (\ref{VB11}) converges uniformly on each compact of $B_{m}$.

For $k\neq\ell$, $f_{k}$ and $f_{\ell}$ are orthogonal in $H\left(
B_{m}\right)  $. This implies that $H\left(  B_{m}\right)  $ is the Hilbert
direct sum%
\[
H\left(  B_{m}\right)  =\widehat{\bigoplus}_{k\geq0}H_{k}\left(  B_{m}\right)
,
\]
where $H_{k}\left(  B_{m}\right)  $ is the space of $k$-homogeneous
polynomials, endowed with the scalar product of $H\left(  B_{m}\right)  $. For
each $k\in\mathbb{N}$, let $\left(  \phi_{k,j}\right)  _{j\in J(k)}$ be an
orthonormal basis of $H_{k}\left(  B_{m}\right)  $; then%
\[
\left(  \phi_{k,j}\right)  _{k\in\mathbb{N},j\in J(k)}%
\]
is a Hilbert basis of $H\left(  B_{m}\right)  $ and
\begin{equation}
\mathcal{K}_{B_{m}}(Z)=\sum_{k\in\mathbb{N},j\in J(k)}\left|  \phi
_{k,j}(Z)\right|  ^{2}=\sum_{k\geq0}\mathcal{K}_{B_{m},k}(Z), \label{VB9}%
\end{equation}
where%
\[
\mathcal{K}_{B_{m},k}(Z)=\sum_{j\in J(k)}\left|  \phi_{k,j}(Z)\right|  ^{2}%
\]
is the reproducing kernel of $H_{k}\left(  B_{m}\right)  $. The expansions
(\ref{VB9}) also converge uniformly on compact subsets of $B_{m}$. Clearly,
$\mathcal{K}_{B_{m},k}$ is a real polynomial, homogeneous of bidegree $\left(
k,k\right)  $. Comparing with the expansion of (\ref{VB8})%
\[
\mathcal{K}_{B_{m}}(Z)=\sum_{k=0}^{\infty}\binom{k+m}{m}\left\|  Z\right\|
^{2k},
\]
we conclude that
\begin{equation}
\mathcal{K}_{B_{m},k}(Z)=\binom{k+m}{m}\left\|  Z\right\|  ^{2k}. \label{VB7}%
\end{equation}

More generally, consider the Hermitian ball $B_{m}(\rho)$ of radius $\rho$.
Its Bergman kernel w.r. to\emph{ the same }$\omega_{m}$ is
\[
\mathcal{K}_{B_{m}(\rho)}(Z)=\frac{1}{\rho^{2m}}\mathcal{K}_{B_{m}}\left(
\frac{Z}{\rho}\right)  ;
\]
the component of bidegree $\left(  k,k\right)  $ is
\begin{equation}
\mathcal{K}_{B_{m}(\rho),k}(Z)=\frac{1}{\rho^{2m+2k}}\binom{k+m}{m}\left\|
Z\right\|  ^{2k}. \label{VB7'}%
\end{equation}
If $f\in H\left(  B_{m}(\rho)\right)  $, its component of degree $k$ is then
given by
\begin{equation}
f_{k}(Z)=\int_{B_{m}(\rho)}\frac{1}{\rho^{2m+2k}}\binom{k+m}{m}\left\langle
Z,W\right\rangle ^{k}f(W)\omega_{V}(W). \label{VB10}%
\end{equation}

\subsection{}

Now we show that the virtual Bergman kernel of $\left(  \Omega,p\right)  $
allows us to compute the Bergman kernel of any inflated domain
\[
\widehat{\Omega}_{m}=\left\{  \left(  z,Z\right)  \in\Omega\times
\mathbb{C}^{m}\mid\left\|  Z\right\|  ^{2}<p(z)\right\}  .
\]
Here $\widehat{\Omega}_{m}$ is endowed with the volume form
\[
\omega_{V}(z)\wedge\omega_{m}(Z),
\]
where $\omega_{m}(Z)=\left(  \frac{\operatorname*{i}}{2\pi}\partial
\overline{\partial}\left\|  Z\right\|  ^{2}\right)  ^{m}$.

\begin{theorem}
\label{TH1}The Bergman kernel $\widehat{K}_{m}$ ($\widehat{\mathcal{K}}_{m}$)
of $\widehat{\Omega}_{m}$ is
\begin{align}
\widehat{K}_{m}\left(  (z,Z),(w,W)\right)   &  =L_{m}\left(  z,w;\left\langle
Z,W\right\rangle \right)  ,\label{VB5'}\\
\widehat{\mathcal{K}}_{m}(z,Z)  &  =\mathcal{L}_{m}\left(  z,\left\Vert
Z\right\Vert ^{2}\right)  ,\quad\label{VB5}%
\end{align}
where
\begin{align}
L_{m}(z,w;r)  &  =\frac{1}{m!}\frac{\partial^{m}}{\partial r^{m}}%
L_{0}(z,w;r),\label{VB12'}\\
\mathcal{L}_{m}(z,r)  &  =\frac{1}{m!}\frac{\partial^{m}}{\partial r^{m}%
}\mathcal{L}_{0}(z,r). \label{VB12}%
\end{align}

\end{theorem}

Note that the relation between $\widehat{\mathcal{K}}_{1}$ and $\widehat
{\mathcal{K}}_{m}$, deduced from (\ref{VB5}) and (\ref{VB12}), is nothing else
that the inflation principle (\ref{VB0}).

\proof We have
\[
L_{m}\left(  z,w;r\right)  =\sum_{k=0}^{\infty}\binom{k+m}{m}K_{\Omega
,p^{k+m}}(z,w)r^{k}.
\]
So we want to prove that the Bergman kernel of $\widehat{\Omega}_{m}$ is
\[
\widehat{K}_{m}\left(  (z,Z),(w,W)\right)  =\sum_{k=0}^{\infty}\binom{k+m}%
{m}K_{\Omega,p^{k+m}}(z,w)\left\langle Z,W\right\rangle ^{k}.
\]

Let $H_{k}\left(  \widehat{\Omega}_{m}\right)  $ be the subspace of functions
$f(z,Z)$ in $H\left(  \widehat{\Omega}_{m}\right)  $, which are $k$%
-homogeneous polynomial w.r.~to the variable $Z$. For $k\neq\ell$, $f\in
H_{k}\left(  \widehat{\Omega}_{m}\right)  $, $g\in H_{\ell}\left(
\widehat{\Omega}_{m}\right)  $, we have
\begin{align*}
\left(  f\mid g\right)  _{\widehat{\Omega}_{m}}  &  =\int_{\widehat{\Omega
}_{m}}f(w,W)\overline{g(w,W)}\omega_{V}(w)\wedge\omega_{m}(W)\\
&  =\int_{w\in\Omega}\left(  f(w,\ )\mid g(w,\ )\right)  _{B_{m}\left(
p(w)^{1/2}\right)  }\omega_{V}(w)=0.
\end{align*}
This implies that $H\left(  \widehat{\Omega}_{m}\right)  $ is the Hilbert
direct sum
\begin{equation}
H\left(  \widehat{\Omega}_{m}\right)  =\widehat{\bigoplus_{k\geq0}}%
H_{k}\left(  \widehat{\Omega}_{m}\right)  . \label{VB16}%
\end{equation}

Fix $k\in\mathbb{N}$. Let $f\in H_{k}\left(  \widehat{\Omega}_{m}\right)  $.
For almost all $w\in\Omega$, the function $W\mapsto f(w,W)$ belongs to
$H\left(  B_{m}\left(  p(w)^{1/2}\right)  \right)  $; by (\ref{VB10}),%
\[
p(w)^{m+k}f(w,Z)=\int_{\left\Vert W\right\Vert ^{2}<p(w)}\binom{k+m}%
{m}\left\langle Z,W\right\rangle ^{k}f(w,W)\omega_{m}(W).
\]
By the reproducing property of $K_{\Omega,p^{k+m}}$, we have%
\[
f\left(  z,Z\right)  =\int_{w\in\Omega}K_{\Omega,p^{k+m}}(z,w)p(w)^{m+k}%
f(w,Z)\omega_{V}(w).
\]
These relations imply%
\[
f(z,Z)=\int_{\widehat{\Omega}_{m}}K_{\Omega,p^{k+m}}(z,w)\binom{k+m}%
{m}\left\langle Z,W\right\rangle ^{k}f(w,W)\omega_{V}(w)\wedge\omega_{m}(W),
\]
which means that
\[
K_{\Omega,p^{k+m}}(z,w)\binom{k+m}{m}\left\langle Z,W\right\rangle ^{k}%
\]
is the reproducing kernel of $H_{k}\left(  \widehat{\Omega}_{m}\right)  $.

From (\ref{VB16}), we deduce that
\[
\sum_{k=0}^{\infty}K_{\Omega,p^{k+m}}(z,w)\binom{k+m}{m}\left\langle
Z,W\right\rangle ^{k}%
\]
is the reproducing kernel of $H\left(  \widehat{\Omega}_{m}\right)  $.
\endproof

\section{Virtual Bergman kernels for bounded symmetric domains}

In this section, we compute the virtual Bergman kernel $\mathcal{L}_{\Omega
,p}\left(  z,r\right)  =\mathcal{L}_{0}\left(  z,r\right)  $ when $\Omega$ is
an irreducible bounded circled symmetric domain and $p$ is a power of the
generic norm of $\Omega$.

\subsection{}

Let $V$ be a complex finite-dimensional vector space and $\Omega\subset V$ an
irreducible bounded circled symmetric domain. Then $V$ is endowed with a
canonical structure of \emph{positive Hermitian Jordan triple}. The numerical
invariants of $V$ (or of $\Omega$) are the \emph{rank} $r$ and the
\emph{multiplicities} $a$ and $b$ ($b=0$ iff the domain is of tube type); the
\emph{genus }is defined by
\[
g=2+a(r-1)+b.
\]
The \emph{generic minimal polynomial} $m\left(  T,x,y\right)  $ and the
\emph{generic norm} $N(x,y)  $ of $\Omega$ are written as
\begin{align*}
m\left(  T,x,y\right)   &  =T^{r}-T^{r-1}m_{1}(x,y)+\cdots+(-1)^{r}%
m_{r}(x,y),\\
N\left(  x,y\right)   &  =m\left(
1,x,y\right)=1-m_{1}(x,y)+\cdots+(-1)^{r}m_{r}(x,y),
\end{align*}
where $m_{1},\ldots,m_{r}$ are polynomials on $V\times\overline{V}$,
homogeneous of bidegrees $\left(  1,1\right)  ,\ldots,$ $\left(  r,r\right)  $
respectively. In particular, $m_{1}$ is an Hermitian inner product on $V$; we
endow $V$ with the K\"{a}hler form
\[
\alpha(z)=\frac{\operatorname{i}}{2\pi}\partial\overline{\partial}m_{1}(z,z)
\]
and with the volume form
\[
\omega=\alpha^{n},
\]
where $n$ is the complex dimension of $V$.

(See \cite{YinRoos2003} for a review of the above properties).

\subsection{}

The \emph{Bergman kernel} of $\Omega$ is then
\begin{equation}
\mathcal{K}(z)=\frac{1}{\operatorname{vol}\Omega}\frac{1}{N\left(  z,z\right)
^{g}},
\end{equation}
with $\operatorname{vol}\Omega=\int_{\Omega}\omega$.

More generally, consider the \emph{weighted Bergman space} of $\Omega$ with
respect to a power of the generic norm:%
\[
H^{(\mu)}(\Omega)=\left\{  f\in\operatorname{Hol}\Omega\mid\int_{\Omega
}\left\vert f(z)\right\vert ^{2}N(z,z)^{\mu}\omega(z)<\infty\right\}  .
\]
For $\mu>-1$, the space $H^{(\mu)}(\Omega)$ is non-zero and is a Hilbert space
of holomorphic functions. Its reproducing kernel is
\begin{equation}
\mathcal{K}^{(\mu)}(z)=\frac{1}{\int_{\Omega}N(z,z)^{\mu}\omega(z)}N\left(
z,z\right)  ^{-g-\mu}. \label{WeightBK}%
\end{equation}

The denominator $\int_{\Omega}N(z,z)^{\mu}\omega(z)$ of the previous formula is
called the \emph{Hua integral}. It has been computed for the four series of
classical domains (with different normalizations of the volume element) by Hua
L.K. \cite{Hua1963}.

\begin{theorem}
\label{TH2}\cite{YinRoos2003}Let $\Omega$ be an irreducible bounded circled
symmetric domain. The value of the Hua integral is given by%
\begin{equation}
\int_{\Omega}N(z,z)^{s}\omega(z)=\frac{\chi(0)}{\chi(s)}\int_{\Omega}\omega,
\label{HuaInt}%
\end{equation}
where $\chi$ is the polynomial of degree $n=\dim_{\mathbb{C}}\Omega$, related
to the numerical invariants of $\Omega$ by
\begin{equation}
\chi(s)=\prod\limits_{j=1}^{r}%
\left(  s+1+(j-1)\frac{a}{2}\right)  _{1+b+(r-j)a}. \label{HuaPol}%
\end{equation}

\end{theorem}

Here $\left(  s\right)  _{k}$ denotes the \emph{raising factorial}%
\[
(s)_{k}=s(s+1)\cdots(s+k-1)=\frac{\Gamma(s+k)}{\Gamma(s)}.
\]

The proof uses the polar decomposition in positive Hermitian Jordan triples
(which generalizes the polar decomposition of matrices) and the following
generalization, due to Selberg \cite{Selberg1944}, of the Beta integral:%
\begin{align}
&  \int_{0}^{1}\cdots\int_{0}^{1}%
{\displaystyle\prod\limits_{j=1}^{n}}
t_{j}^{x-1}\left(  1-t_{j}\right)  ^{y-1}%
{\displaystyle\prod\limits_{1\leq j<k\leq n}} \left\vert
t_{j}-t_{k}\right\vert ^{2z}\operatorname{d}t_{1}\cdots
\operatorname{d}t_{n}\nonumber\\
&  =%
{\displaystyle\prod\limits_{j=1}^{n}}
\frac{\Gamma(x+(j-1)z)\Gamma(y+(j-1)z)\Gamma(jz+1)}{\Gamma(x+y+(n+j-2)z)\Gamma
(z+1)},
\end{align}
for $\operatorname{Re}x>0$, $\operatorname{Re}y>0$, $\operatorname{Re}%
z>\min\left(  \frac{1}{n},\frac{\operatorname{Re}x}{n-1},\frac
{\operatorname{Re}y}{n-1}\right)  $.

\begin{corollary}
The Bergman kernel of $H^{(\mu)}(\Omega)$ is
\begin{equation}
\mathcal{K}^{(\mu)}(z)=\frac{\chi(\mu)}{\chi(0)}N\left(  z,z\right)  ^{-\mu
}\mathcal{K}(z), \label{WeightBK2}%
\end{equation}
where $\mathcal{K}=\mathcal{K}^{(0)}$ is the Bergman kernel of $\Omega$ .
\end{corollary}

\subsection{}

\begin{theorem}
\label{TH3}Let $\Omega$ be an irreducible bounded circled symmetric domain.
The virtual Bergman kernel of $\left(  \Omega,N(z,z)^{\mu}\right)  $ is
\begin{equation}
\mathcal{L}^{(\mu)}\left(  z,r\right)  =\mathcal{K}(z)F_{\chi,\mu}\left(
\frac{r}{N(z,z)^{\mu}}\right)  , \label{VBKern}%
\end{equation}
where $\mathcal{K}$ is the Bergman kernel of $\Omega$, $\chi$ the polynomial
defined by (\ref{HuaPol}) and $F_{\chi,\mu}$ is the rational function
\begin{equation}
F_{\chi,\mu}(t)=\frac{1}{\chi(0)}%
{\displaystyle\sum\limits_{k=0}^{\infty}}
\chi(\mu k)t^{k}. \label{VBRat}%
\end{equation}

\end{theorem}

The proof of (\ref{VBKern}) is straightforward, using (\ref{VB3}) and
(\ref{WeightBK2}):%
\begin{align*}
\mathcal{L}^{(\mu )}\left(  z,r\right)   &  =\sum_{k=0}^{\infty}%
\mathcal{K}^{(k\mu)}(z)r^{k}\\
&  =\sum_{k=0}^{\infty}\frac{\chi(k\mu)}{\chi(0)}N\left(  z,z\right)  ^{-k\mu
}\mathcal{K}(z)r^{k}\\
&  =\frac{\mathcal{K}(z)}{\chi(0)}\sum_{k=0}^{\infty}\chi(k\mu)\left(
\frac{r}{N(z,z)^{\mu}}\right)  ^{k}.
\end{align*}

If the polynomial $k\mapsto\chi(k\mu)$ is decomposed as
\begin{equation}
\frac{\chi(k\mu)}{\chi(0)}=\sum_{j=0}^{n}c_{\mu,j}\frac{\left(  k+1\right)
_{j}}{j!}, \label{VBComb}%
\end{equation}
the function defined by (\ref{VBRat}) is
\[
F_{\chi,\mu}(t)=\sum_{j=0}^{n}c_{\mu,j}\left(  \frac{1}{1-t}\right)  ^{j}.
\]

\section{Tables for bounded symmetric domains}

Hereunder we give the results for each type of irreducible bounded symmetric
domains, which allow the reader to apply Theorems \ref{TH1} and \ref{TH3} to
special cases.

\subsection{Classification of irreducible circled bounded symmetric domains}

Here is the complete list of irreducible circled bounded symmetric domains, up
to linear isomorphisms. There is some overlapping between the four infinite
families, due to a finite number of isomorphisms in low dimemsions.

\subsubsection*{\textit{Type I}$_{m,n}$\textit{ }$\left(  1\leq m\leq
n\right)  $}

$V=\mathcal{M}_{m,n}(\mathbb{C})$ (space of $m\times n$ matrices with complex
entries).
\[
\Omega=\left\{  x\in V\mid I_{m}-x^{t}\overline{x}\gg0\right\}  \mathit{.}%
\]

\subsubsection*{\textit{Type II}$_{n}$\textit{ }$\left(  n\geq2\right)  $}

$V=\mathcal{A}_{n}(\mathbb{C})$ (space of $n\times n$ alternating matrices).
\[
\Omega=\left\{  x\in V\mid I_{n}+x\overline{x}\gg0\right\}  .
\]

\subsubsection*{\textit{Type III}$_{n}$\textit{ }$\left(  n\geq1\right)  $}

$V=\mathcal{S}_{n}(\mathbb{C})$ (space of $n\times n$ symmetric matrices).
\[
\Omega=\left\{  x\in V\mid I_{n}-x\overline{x}\gg0\right\}  .
\]

\subsubsection*{\textit{Type IV}$_{n}$\textit{ }$\left(  n\neq2\right)  $}

$V=\mathbb{C}^{n}$, $q(x)=\sum x_{i}^{2},$ $q(x,y)=2\sum x_{i}y_{i}$. The
domain $\Omega$ is defined by
\[
1-q(x,\overline{x})+\left\vert q(x)\right\vert ^{2}>0,\quad2-q(x,\overline
{x})>0.
\]

\subsubsection*{\textit{Type V}}

$V=\mathcal{M}_{2,1}(\mathbb{O}_{\mathbb{C}})\simeq\mathbb{C}^{16}$,
exceptional type.

\subsubsection*{\textit{Type VI}}

$V=\mathcal{H}_{3}(\mathbb{O}_{\mathbb{C}})\simeq\mathbb{C}^{27}$, exceptional type.

\subsection{Numerical invariants}

\subsubsection*{\textit{Type I}$_{m,n}$\textit{ }$\left(  1\leq m\leq
n\right)  $}%

\[
r=m,\quad a=2,\quad b=n-m,\quad g=m+n.
\]

\subsubsection*{\textit{Type II}$_{2p}$ $(p\geq1)$}%

\[
r=\frac{n}{2}=p,\quad a=4,\quad b=0,\quad g=2\left(  n-1\right)  .
\]

\subsubsection*{\textit{Type II}$_{2p+1}$ $(p\geq1)$}%

\[
r=\left[  \frac{n}{2}\right]  =p,\quad a=4,\quad b=2,\quad g=2(n-1).
\]

\subsubsection*{\textit{Type III}$_{n}$\textit{ }$\left(  n\geq1\right)  $}%

\[
r=n,\quad a=1,\quad b=0,\quad g=n+1.
\]

\subsubsection*{\textit{Type IV}$_{n}$\textit{ }$\left(  n\neq2\right)  $}%

\[
r=2,\quad a=n-2,\quad b=0,\quad g=n.
\]

\subsubsection*{\textit{Type V}}%

\[
r=2,\quad a=6,\quad b=4,\quad g=12.
\]

\subsubsection*{\textit{Type VI}}%

\[
r=3,\quad a=8,\quad b=0,\quad g=18.
\]

\subsection{Generic norm}

\subsubsection*{\textit{Type I}$_{m,n}$\textit{ }$\left(  1\leq m\leq
n\right)  $}

$V=\mathcal{M}_{m,n}(\mathbb{C})$ (space of $m\times n$ matrices with complex
entries).
\[
N(x,y)=\operatorname*{Det}(I_{m}-x^{t}\overline{y}).
\]

\subsubsection*{\textit{Type II}$_{n}$\textit{ }$\left(  n\geq2\right)  $}

$V=\mathcal{A}_{2p}(\mathbb{C})$ (space of $n\times n$ alternating matrices).
\[
N(x,y)^{2}=\operatorname*{Det}(I_{n}+x\overline{y}).
\]

\subsubsection*{\textit{Type III}$_{n}$\textit{ }$\left(  n\geq1\right)  $}

$V=\mathcal{S}_{n}(\mathbb{C})$ (space of $n\times n$ symmetric matrices).
\[
N(x,y)=\operatorname*{Det}(I_{n}-x\overline{y}).
\]

\subsubsection*{\textit{Type IV}$_{n}$\textit{ }$\left(  n\neq2\right)  $}

$V=\mathbb{C}^{n}$.
\[
N(x,y)=1-q(x,\overline{y})+q(x)q(\overline{y}).
\]

\subsubsection*{\textit{Type V}}

$V=\mathcal{M}_{2,1}(\mathbb{O}_{\mathbb{C}})$.%
\[
N(x,y)=1-(x|y)+(x^{\sharp}|y^{\sharp}).
\]

\subsubsection*{\textit{Type VI}}

$V=\mathcal{H}_{3}(\mathbb{O}_{\mathbb{C}})$.
\[
N(x,y)=1-(x|y)+(x^{\sharp}|y^{\sharp})-\det x\det\overline{y}.
\]

\subsection{The polynomial $\chi$ for the Hua integral\label{HuaPolList}}

Recall that for $s>-1$,%

\[
\chi(s)\int_{\Omega}N(x,x)^{s}\omega=\chi(0)\int_{\Omega}\omega.
\]

\subsubsection*{\textit{Type I}$_{m,n}$}%

\[
{\chi(s)=\prod_{j=1}^{m}(s+j)_{n}.}%
\]

\subsubsection*{\textit{Type II}$_{2p}$}%

\[
{\chi(s)=\prod_{j=1}^{p}(s+2j-1)_{2p-1}.}%
\]

\subsubsection*{\textit{Type II}$_{2p+1}$}%

\[
{\chi(s)=\prod_{j=1}^{p}(s+2j-1)_{2p+1}.}%
\]

\subsubsection*{\textit{Type III}$_{n}$}%

\[
{\chi(s)=\prod_{j=1}^{n}\left(  s+\frac{j+1}{2}\right)  _{1+n-j}.}%
\]

\subsubsection*{\textit{Type IV}$_{n}$}%

\[
{\chi(s)=\left(  s+1\right)  _{n-1}\left(  s+\frac{n}{2}\right)  .}%
\]

\subsubsection*{\textit{Type V}}%

\[
{\chi(s)=(s+1)_{8}(s+4)_{8}.}%
\]

\subsubsection*{\textit{Type VI}}%

\[
{\chi(s)=(s+1)_{9}(s+5)_{9}(s+9)_{9}.}%
\]

\section{Open problems}

\subsection{}

Understand the rational function%
\[
F_{\chi,\mu}(t)=\frac{1}{\chi(0)}%
%TCIMACRO{\dsum \limits_{k=0}^{\infty}}%
%BeginExpansion
{\displaystyle\sum\limits_{k=0}^{\infty}}
%EndExpansion
\chi(\mu k)t^{k}=\sum_{j=0}^{n}c_{\mu,j}\left(  \frac{1}{1-t}\right)  ^{j},
\]
when $\chi$ is a polynomial from the list of subsection \ref{HuaPolList}.
Recall that the coefficients $c_{\mu,j}$ are given by
\[
\frac{\chi(k\mu)}{\chi(0)}=\sum_{j=0}^{n}c_{\mu,j}\frac{\left(  k+1\right)
_{j}}{j!}.
\]

\subsection{}

The virtual Bergman kernel
\[
\mathcal{L}_{\Omega,p}\left(  z,r\right)  =\sum_{k=0}^{\infty}\mathcal{K}%
_{\Omega,p^{k}}(z)r^{k}%
\]
is suitable for the computation of the Bergman kernel of \textquotedblleft
inflated domains by Hermitian balls\textquotedblright%
\[
\widehat{\Omega}_{m}=\left\{  \left(  z,Z\right)  \in\Omega\times
\mathbb{C}^{m}\mid\left\Vert Z\right\Vert ^{2}<p(z)\right\}  ,
\]
which may also be written%
\[
\widehat{\Omega}_{m}=\left\{  \left(  z,Z\right)  \in\Omega\times
\mathbb{C}^{m}\mid Z\in\left(  p(z)\right)  ^{1/2}B_{m}\right\}  .
\]

If $F$ is any circled domain in $\mathbb{C}^{m}$, what can be said about the
Bergman kernel of
\[
\widehat{\Omega}_{F}=\left\{  \left(  z,Z\right)  \in\Omega\times
\mathbb{C}^{m}\mid Z\in\left(  p(z)\right)  ^{1/2}F\right\}  ,
\]
for example when $\Omega$ is a bounded symmetric domain, $p=N_{\Omega
}(z,z)^{\mu}$ and $F$ another bounded symmetric domain?

For some families $\left\{  F\right\}  $ other than the family $\left\{
B_{m}\right\}  $ of Hermitian balls, is it possible to define an analogous of
the virtual Bergman kernel $\mathcal{L}_{\Omega,p}$ and obtain an analogous of
Theorem \ref{TH1}?

\subsection{}

It would also be interesting to replace the weighted Bergman space
$H(\Omega,p)$ by more general Hilbert spaces of holomorphic functions; for
example, when $\Omega$ is a bounded symmetric domain, the spaces with
reproducing kernel $N(z,z)^{\mu}$, where $\mu$ is in the \textquotedblright
Berezin-Wallach set\textquotedblright\ of $\Omega$ (see
\cite{FarautKoranyi1990}).

\end{document}